\documentclass[11pt]{article}

\usepackage[T1]{fontenc}
\usepackage[utf8]{inputenc}
\usepackage[margin=1in]{geometry}
\usepackage{amsmath,amssymb,amsthm,mathtools}
\usepackage{enumitem}
\usepackage{microtype}
\usepackage{hyperref}

\hypersetup{
  colorlinks=true,
  linkcolor=blue,
  citecolor=blue,
  urlcolor=blue
}

\numberwithin{equation}{section}

\newtheorem{theorem}{Theorem}[section]
\newtheorem{lemma}[theorem]{Lemma}
\newtheorem{proposition}[theorem]{Proposition}

\theoremstyle{definition}

\newtheorem{example}[theorem]{Example}

\theoremstyle{remark}

\newcommand{\Z}{\mathbb Z}
\newcommand{\IAM}{\operatorname{IAM}}
\newcommand{\wt}{\operatorname{wt}}
\newcommand{\dist}{\operatorname{dist}}

\title{
The Entire Integer-antimagic\\ Spectra
of Connected Graphs
}

\author{%
  Richard M. Low\thanks{Department of Mathematics and Statistics, San Jose State University, San Jose, CA 95192, USA.
    Email: \href{mailto:richard.low@sjsu.edu}{\texttt{richard.low@sjsu.edu}}.}
  \and
  U\u{g}ur Odaba\c{s}i\thanks{Department of Engineering Sciences, Istanbul University-Cerrahpasa, Istanbul 34320, Turkiye.
    Email: \href{mailto:ugur.odabasi@iuc.edu.tr}{\texttt{ugur.odabasi@iuc.edu.tr}}.}
  \and
  Dan Roberts\thanks{Department of Mathematics, Illinois Wesleyan University, Bloomington, IL 61701, USA.
    Email: \href{mailto:drobert1@iwu.edu}{\texttt{drobert1@iwu.edu}}.}
}
\date{November 26, 2026}

\begin{document}
\maketitle

\begin{abstract}
For a connected simple graph $G$, let $\IAM(G)$
denote the set of integers $k\geq2$ for which
$G$ admits a nowhere-zero $\Z_k$-antimagic
edge labeling.

We prove that every connected simple graph
of order at least three has an integer-antimagic
spectrum consisting of all integers above its
least admissible modulus. Combined with
the known minimum-modulus classification of
Jones and Zhang, this gives the complete
integer-antimagic spectrum of every such graph
and establishes the interval-spectrum conjecture.

The principal new ingredient is a constructive
theorem for trees: every tree of order $n\geq3$
is $\Z_k$-antimagic for $k\geq n+2$,
and the bound improves to $k\geq n+1$
except for even-order stars.

The proof uses feasible vertex weights,
path switching, and a zero-edge elimination
method. When a unique zero vertex weight
occurs at an internal vertex, a matching
argument first repairs its incident edges.
The edge-addition theorem of Odabasi,
Roberts, and Low then transfers the tree
result to all connected graphs.
\end{abstract}



\section{Introduction and main results}

Throughout, graphs are finite, simple, and
undirected. For an integer $k\geq2$, let
$\Z_k=\Z/k\Z$ denote the additive cyclic
group of order $k$.

A \emph{nowhere-zero $\Z_k$-antimagic
edge labeling} of a graph $G=(V,E)$ is a
function
\[
f:E\longrightarrow\Z_k\setminus\{0\}
\]
whose induced vertex-weight function
\[
\wt_f(v)=\sum_{e\ni v}f(e)
\]
is injective. Edge labels themselves need
not be distinct.

The \emph{integer-antimagic spectrum} is
\[
\IAM(G)=
\{k\geq2:G\text{ is }\Z_k\text{-antimagic}\}.
\]

Since $|V|=n$ distinct vertex weights must
belong to $\Z_k$, every admissible modulus
satisfies $k\geq n$.

For a connected graph $G$ of order at
least three, write
\[
\mu(G)=\min\IAM(G).
\]

Jones and Zhang \cite{JonesZhang}
determined this minimum for every
connected graph of order at least three.
In the terminology used here, their
classification is the following.

\begin{theorem}[Jones--Zhang
\cite{JonesZhang}]
\label{thm:knownminimum}
Let $G$ be a connected simple graph
of order $n\geq3$. Then
\begin{equation}
\mu(G)=
\begin{cases}
n+2,
 &G\cong K_{1,n-1},\ n\text{ even},\\
n+1,
 &G\cong K_3,\\
n+1,
 &n\equiv2\pmod4,\
  G\not\cong K_{1,n-1},\\
n,
 &\text{otherwise}.
\end{cases}
\label{eq:minclassification}
\end{equation}
\end{theorem}

The cases in
\eqref{eq:minclassification}
are read in order; the first case
takes precedence over the third.

Determining $\mu(G)$ does not by
itself determine $\IAM(G)$.
Admissibility at a modulus $k$
need not imply admissibility at
$k+1$.

The conjecture that these spectra
have no gaps above their least
element appears in earlier work
and is explicitly discussed by
Odabasi, Roberts, and Low
\cite{ORL}.

Our main result establishes this
property for every connected
simple graph of order at least three.

\begin{theorem}[Complete spectra]
\label{thm:complete}
Let $G$ be a connected simple graph
of order $n\geq3$. Then
\begin{equation}
\IAM(G)=
\{k\in\Z:k\geq\mu(G)\},
\label{eq:complete}
\end{equation}
where $\mu(G)$ is given by
\eqref{eq:minclassification}.
\end{theorem}

Thus the entire spectrum depends
only on the order of the graph
and the exceptional conditions
in Theorem~\ref{thm:knownminimum}.

The new construction underlying
Theorem~\ref{thm:complete} is a
sharp sufficient bound for trees.

\begin{theorem}[Sharp tree bounds]
\label{thm:trees}
Let $T$ be a tree of order $n\geq3$.
Then
\begin{equation}
[n+2,\infty)\cap\Z
\subseteq\IAM(T).
\label{eq:alltrees}
\end{equation}

If $n$ is odd or $T$ is not a star,
then
\begin{equation}
[n+1,\infty)\cap\Z
\subseteq\IAM(T).
\label{eq:nonstars}
\end{equation}

Both thresholds are optimal as
uniform sufficient bounds for
the corresponding classes.
\end{theorem}

To pass from trees to connected
graphs, we use the edge-addition
theorem of Odabasi, Roberts, and
Low \cite{ORL}. We record the
connected case needed here.

\begin{theorem}[Edge-addition preservation
\cite{ORL}]
\label{thm:knownstability}
Let $H$ be a connected simple graph
of order at least four, and let
$H'$ be obtained by adding an edge
between two nonadjacent vertices.
Then
\[
\IAM(H)\subseteq\IAM(H').
\]
\end{theorem}

The order-three exception in
the general edge-addition result
is handled separately in our proof.

The proof of
Theorem~\ref{thm:trees} is
self-contained and constructive.
It occupies Sections 2--4.
Section 5 deduces
Theorem~\ref{thm:complete} and
establishes sharpness.
Section 6 gives examples
and concluding remarks.

The published theorems quoted
above are used as external inputs;
their original lengthy proofs
are not reproduced here.

\section{Feasible weights and switching}

Let $T=(V,E)$ be a tree and fix
its bipartition
\[
V=X\sqcup Y.
\]

A prescribed vertex weighting
$w:V\to\Z_k$ is called
\emph{feasible} if
\begin{equation}
\sum_{x\in X}w(x)=
\sum_{y\in Y}w(y).
\label{eq:feasible}
\end{equation}

Feasibility is the only
obstruction to inducing
prescribed vertex weights
by an edge labeling of a tree.

\begin{lemma}[Inverse incidence]
\label{lem:incidence}
A weighting $w:V\to\Z_k$
is induced by an edge labeling
$f:E\to\Z_k$, temporarily allowing
zero labels, if and only if $w$
is feasible.

The inducing edge labeling,
when it exists, is unique.
\end{lemma}

\begin{proof}
Every edge has one endpoint
in each bipartition class.
Hence an induced weighting
satisfies
\[
\sum_{x\in X}w(x)
=\sum_{e\in E}f(e)
=\sum_{y\in Y}w(y).
\]

Conversely, let $u$ be a leaf
with neighbor $v$.
Its incident edge must have
label $w(u)$.

Delete $u$ and replace the
prescribed weight at $v$
by $w(v)-w(u)$.

Repeat this leaf-elimination
procedure. The signed
difference of the two
bipartition sums remains zero.

Thus the final remaining
vertex has prescribed
residual weight zero,
and the labels obtained
induce the original weighting.

Every edge label was forced
at some stage, proving
uniqueness.
\end{proof}

We shall repeatedly interchange
weights of vertices belonging to
the same bipartition class.
Such interchanges preserve
feasibility and the multiset
of vertex weights.

The effect on the induced
edge labels is especially simple.

\begin{lemma}[Path switching]
\label{lem:switch}
Let $w$ be feasible and let
$f$ be its inducing edge labeling.

Suppose $p$ and $q$ belong
to the same bipartition class.
Interchange $w(p)$ and $w(q)$,
and put
\[
\delta=w(q)-w(p).
\]

Let
\[
p=x_0,x_1,\ldots,x_{2t}=q
\]
be the unique path between them.

The new inducing labeling $f'$
satisfies
\begin{equation}
f'(x_{j-1}x_j)
=
f(x_{j-1}x_j)
+(-1)^{j-1}\delta
\quad (1\leq j\leq2t),
\label{eq:switch}
\end{equation}
and every edge outside the
path retains its label.

If $w$ is injective, then
$\delta\neq0$.
\end{lemma}

\begin{proof}
The path has even length
because its endpoints lie
in the same bipartition class.

Modify the labels along it
by the alternating increments
in \eqref{eq:switch}.

The induced weight at $p$
increases by $\delta$, while
that at $q$ decreases by
$\delta$.

At each internal vertex,
the changes on its two
path edges cancel.

Hence the modified edge
labels induce precisely
the swapped weighting.
Uniqueness follows from
Lemma~\ref{lem:incidence}.

The final assertion follows
from injectivity.
\end{proof}

In particular, if $p-v-q$
is a path of length two,
then
\begin{equation}
f'(pv)=f(pv)+\delta,\qquad
f'(vq)=f(vq)-\delta,
\label{eq:twoswitch}
\end{equation}
and no other edge label changes.

For a path of length four,
the changes are
\[
+\delta,\quad-\delta,\quad
+\delta,\quad-\delta.
\]

These two switching identities
are the only transformations
needed in the proof.

\section{Eliminating zero edge labels}

Prescribing feasible weights
does not guarantee that the
inducing edge labels are
all nonzero.

The next theorem shows that
zero labels can be eliminated
under conditions satisfied
by our constructions.

\begin{theorem}[Zero-edge elimination]
\label{thm:elimination}
Let $T$ be a tree and let
$w:V(T)\to\Z_k$ be an injective
feasible weighting.

Suppose either
\begin{enumerate}[
label=\textup{(\roman*)},
leftmargin=2.2em]
\item every vertex weight
is nonzero; or
\item exactly one vertex $r$
has weight zero and
$\deg(r)\geq2$.
\end{enumerate}

Then $T$ admits a
nowhere-zero $\Z_k$-antimagic
edge labeling.

The final vertex weights
are permutations of the
original weights within
their bipartition classes.
In case (ii), the zero
weight remains at $r$.
\end{theorem}

The nonzero-weight case uses
a simple transport operation.
The single-zero case also
requires repairing the edges
incident with the zero-weight
vertex.

\subsection{A matching lemma}

\begin{lemma}[One forbidden value]
\label{lem:matching}
Let $d\geq2$ distinct elements
be assigned bijectively to
$d$ positions, where each
position forbids at most
one of the elements.

A permissible assignment
exists unless all positions
forbid the same element.

In that exceptional situation,
no permissible assignment exists.
\end{lemma}

\begin{proof}
Form the bipartite graph
of allowed assignments:
the left vertices are
the positions, the right
vertices are the elements,
and adjacency means that
an assignment is allowed.

Every position has at
least $d-1$ neighbors.

Hence every nonempty
collection of at most
$d-1$ positions has at
least as many available
elements as positions.

For all $d$ positions,
the union of allowed
elements fails to contain
all $d$ elements precisely
when one element is
forbidden everywhere.

Because each position
forbids at most one
element, this is exactly
the stated exception.

Hall's marriage theorem
now gives the result.
\end{proof}

\subsection{Repair at a zero-weight root}

\begin{lemma}[Root-edge repair]
\label{lem:root}
Let $w$ be an injective
feasible weighting having
exactly one zero weight,
at a vertex $r$ with
$d=\deg(r)\geq2$.

There is a sequence of
same-class weight swaps
that fixes $w(r)=0$,
preserves injectivity and
feasibility, and makes
every edge incident with
$r$ nonzero.
\end{lemma}

\begin{proof}
Root $T$ at $r$ and
write its neighbors as
$v_1,\ldots,v_d$.

Let
\[
a_i=w(v_i),\qquad
\ell_i=f(rv_i),\qquad
c_i=\ell_i-a_i.
\]

The weights $a_i$
are distinct and nonzero.

The constant $c_i$
is the portion of the
root-edge label not
determined by the weight
currently assigned to $v_i$.
It remains fixed when
only the neighbor weights
are permuted.

Indeed, a permutation
assigning $a_{\pi(i)}$
to $v_i$ changes the
root-edge label to
\begin{equation}
\ell_i'=c_i+a_{\pi(i)}.
\label{eq:rootlabels}
\end{equation}

This follows by expressing
the permutation as
transpositions between
neighbors of $r$ and
using \eqref{eq:twoswitch}.

We need a permutation
for which
\[
a_{\pi(i)}\neq-c_i
\qquad(1\leq i\leq d).
\]

Thus position $i$ forbids
at most one available
neighbor weight.

By Lemma~\ref{lem:matching},
the required permutation
exists unless
\begin{equation}
c_1=\cdots=c_d=-a
\label{eq:exception}
\end{equation}
for some $a$ among
$a_1,\ldots,a_d$.

We eliminate this
exceptional configuration.

Since $w(r)=0$,
\[
\sum_{i=1}^{d}\ell_i=0,
\]
and hence
\begin{equation}
\sum_{i=1}^{d}c_i
=
-\sum_{i=1}^{d}a_i.
\label{eq:rootidentity}
\end{equation}

If $d=2$ and
\eqref{eq:exception}
holds, then
\[
a_1+a_2=2a.
\]

As $a$ equals one of
$a_1,a_2$, additive
cancellation implies
$a_1=a_2$, contrary
to injectivity.

Now suppose $d\geq3$
and \eqref{eq:exception}
holds.

No neighbor $v_i$
can be a leaf.
For a leaf neighbor,
\[
\ell_i=w(v_i)=a_i,
\]
which gives $c_i=0$.
Equation~\eqref{eq:exception}
would then force $a=0$,
contrary to the fact that
$a$ is a nonzero neighbor
weight.

Consequently every
neighbor $v_i$ has
a child in the tree
rooted at $r$.

Choose two distinct
branches, with neighbors
$v_i,v_j$, and choose
children $u$ of $v_i$
and $z$ of $v_j$.

Their connecting path is
\[
u-v_i-r-v_j-z.
\]

Interchange $w(u)$
and $w(z)$, and put
\[
\delta=w(z)-w(u)\neq0.
\]

Both vertices belong
to the same bipartition
class as $r$, but
neither equals $r$.

By Lemma~\ref{lem:switch},
the changes along this
four-edge path are
\[
+\delta,-\delta,
+\delta,-\delta.
\]

In particular, the
root-edge labels on
$rv_i$ and $rv_j$
change by $-\delta$
and $+\delta$,
respectively.

The neighbor weights
$a_i$ do not change.
Therefore the constants
become
\[
c_i'=c_i-\delta,\qquad
c_j'=c_j+\delta,
\]
while every other
constant is unchanged.

Choose a third branch
$h\notin\{i,j\}$.
Initially all constants
equaled $-a$, but now
\[
c_i'=-a-\delta
\neq-a=c_h'.
\]

Thus the constants
are no longer all
equal, and the
exceptional matching
configuration is gone.

Apply
Lemma~\ref{lem:matching}
and permute the neighbor
weights accordingly.

Equation~\eqref{eq:rootlabels}
then makes every root
edge nonzero.

The cross-branch swap
may have changed labels
below $r$, possibly
creating zero edges there.
This causes no difficulty:
the next step eliminates
all such zeros.

The zero vertex weight
remains at $r$, and
injectivity and
feasibility are preserved.
\end{proof}

\subsection{Proof of zero-edge elimination}

\begin{proof}[Proof of
Theorem~\ref{thm:elimination}]
Let $f$ be the unique
inducing labeling of $w$.

\medskip
\noindent
\textbf{Case (i): All vertex
weights are nonzero.}

Choose any zero edge $uv$.
Neither endpoint can be
a leaf, since a leaf's
vertex weight equals
its unique incident
edge label.

Choose a neighbor
$x\neq u$ of $v$.
Swap $w(u)$ and $w(x)$,
and write
\[
\delta=w(x)-w(u)\neq0.
\]

By \eqref{eq:twoswitch},
\[
f'(uv)=\delta,\qquad
f'(vx)=f(vx)-\delta.
\]

The zero on $uv$
disappears, and at
most one new zero
appears, on $vx$.

If that new label is
zero, repeat along
$v-x-y$, where $y$
is a neighbor of $x$
other than $v$.

Such a neighbor exists:
otherwise $x$ would
be a leaf incident
with a zero edge,
contrary to the
nonzero vertex weights.

Whenever the zero
survives, it advances
along a non-backtracking
walk in the tree.

A non-backtracking walk
in a tree cannot revisit
a vertex and hence
cannot continue forever.

Therefore the
transported zero
eventually disappears.

Each completed
transport strictly
decreases the total
number of zero edges.
Repeating removes
all zero edges.

\medskip
\noindent
\textbf{Case (ii): One zero
weight at an internal vertex.}

Root $T$ at the unique
zero-weight vertex $r$.

Apply
Lemma~\ref{lem:root}
to make all edges
incident with $r$
nonzero.

Now consider a
remaining zero edge
$uv$, with $u$ the
parent of $v$.

Since all root edges
are nonzero, neither
endpoint is $r$.

The vertex $v$
has nonzero weight,
so it cannot be
a leaf.

Choose a child
$x$ of $v$.
Swap the weights
at $u$ and $x$.

Both vertices
are distinct from $r$
and belong to the
same bipartition class.

With
\[
\delta=w(x)-w(u)\neq0,
\]
we again obtain
\[
f'(uv)=\delta,\qquad
f'(vx)=f(vx)-\delta.
\]

The original zero
disappears.
If a zero survives,
it moves from $uv$
to $vx$.

Let
\[
h(y)=\dist_T(r,y)
\]
be the depth
of a vertex $y$.

The depth of the
child endpoint of
a surviving zero
edge increases
by one at each
transport step.

This depth is bounded
by the height of
the finite tree.

Moreover, a surviving
zero cannot reach
a leaf, because
every vertex other
than $r$ retains
a nonzero weight.

Hence each outward
transport terminates
by eliminating a zero.

All affected edges
lie strictly below
the root. Thus
the previously
repaired root-edge
labels are unchanged.

Repeating the
procedure removes
all remaining zeros.

Throughout both
cases, every operation
merely swaps weights
within a bipartition
class. Thus
feasibility and
injectivity are
preserved.

In the second case,
the unique zero
vertex weight is
never moved.

The final inducing
edge labeling is
nowhere-zero, and
its induced vertex
weights are pairwise
distinct.
\end{proof}

\section{Explicit constructions for trees}

An \emph{ordinary
opposite pair} modulo
$k$ is a pair
\[
\{a,-a\}
\]
whose elements are
distinct and nonzero.

Its two weights sum
to zero.

The number of
ordinary opposite
pairs in $\Z_k$ is
\begin{equation}
P(k)=
\begin{cases}
(k-1)/2,&k\text{ odd},\\
(k-2)/2,&k\text{ even}.
\end{cases}
\label{eq:pairs}
\end{equation}

For even $k$,
the nonzero
self-inverse
residue $k/2$
belongs to no
ordinary opposite
pair.

We use complete
opposite pairs
within bipartition
classes so that
the class sums
remain unchanged.

\begin{lemma}[Feasible vertex weights]
\label{lem:weights}
Let $T$ be a tree of
order $n\geq3$ with
bipartition
$V=X\sqcup Y$.

\begin{enumerate}[
label=\textup{(\roman*)},
leftmargin=2.2em]

\item If $n$ is even
and $T$ is not a star,
then for every
$k\geq n+1$
there is an injective
feasible weighting
with all weights
nonzero.

\item If $n$ is odd
and $k=n+1$, there
is an injective
feasible weighting
with exactly one
zero weight, at
an internal vertex.

\item If $n$ is odd
and $k\geq n+2$,
there is an injective
feasible weighting
with all weights
nonzero.
\end{enumerate}
\end{lemma}

\begin{proof}
We use the five
constructions in
Table~\ref{tab:cases}.
Their feasibility
and residue counts
are verified below.

\begin{table}[ht]
\centering
\small
\begin{tabular}{p{0.31\linewidth}
                p{0.32\linewidth}
                p{0.27\linewidth}}
\hline
Case & Initial assignment & Completion\\
\hline
Even $n$; both
classes even
& None
& Complete opposite pairs\\

Even $n$; both
classes odd
& $(1,2,-3)$ in one
class; $(-1,-2,3)$
in the other
& Unused opposite pairs\\

Odd $n$;
$k=n+1$
& $0$ at an internal
vertex of the
odd-sized class
& All ordinary pairs,
omitting $k/2$\\

Odd $n$;
odd $k\geq n+2$
& $-1$ in the odd
class; $(1,-2)$
in the even class
& Unused opposite pairs\\

Odd $n$;
even $k\geq n+2$
& $k/2$ in the odd
class; $(1,k/2-1)$
in the even class
& Unused opposite pairs\\
\hline
\end{tabular}
\caption{Feasible-weight constructions.
Each initial assignment gives equal partial
sums in the two bipartition classes.}
\label{tab:cases}
\end{table}

\medskip
\noindent
\textbf{(i) Even order,
nonstar tree.}

Let
\[
r=|X|,\qquad
s=|Y|,\qquad
r+s=n.
\]

The two class
sizes have the
same parity.

If both are even,
fill each class
with complete
ordinary opposite
pairs.

Exactly $n/2$
pairs are required.

For odd $k\geq n+1$,
\[
P(k)=\frac{k-1}{2}
\geq\frac n2.
\]

For even $k\geq n+1$,
parity gives
$k\geq n+2$, so
\[
P(k)=\frac{k-2}{2}
\geq\frac n2.
\]

Suppose both
class sizes are
odd.

Because $T$
is not a star,
neither class
has size one.
Thus
\[
r,s\geq3,
\qquad n\geq6.
\]

Assign
\[
1,\quad2,\quad-3
\]
to three vertices
of $X$, and
\[
-1,\quad-2,\quad3
\]
to three vertices
of $Y$.

These six residues
are distinct and
nonzero because
$k\geq n+1\geq7$.

Both triples
sum to zero.

Complete the
remaining vertices
with ordinary
opposite pairs,
avoiding the
three classes
$\{\pm1\}$,
$\{\pm2\}$,
and $\{\pm3\}$.

The number of
additional pairs
is $(n-6)/2$.

Including the
three excluded
pair classes,
the total needed
is
\[
3+\frac{n-6}{2}
=\frac n2
\leq P(k).
\]

Thus every
even-order nonstar
tree has the
required nonzero
feasible weighting.

\medskip
\noindent
\textbf{(ii) Odd order,
$k=n+1$.}

Let $X$ be the
odd-sized
bipartition class
and $Y$ the
even-sized class.

We first show
that $X$ contains
an internal vertex.

If every vertex
of $X$ were a
leaf, a path
between two
distinct vertices
of $Y$ would
contain an
internal vertex
of $X$.

Hence $Y$
would have
at most one
vertex, impossible
because $Y$ is
nonempty and
has even size.

Choose an internal
vertex $r\in X$
and assign
\[
w(r)=0.
\]

Now $k=n+1$
is even.
The residues
other than
$0$ and $k/2$
form exactly
\[
\frac{k-2}{2}
=\frac{n-1}{2}
\]
ordinary
opposite pairs.

Both $|X|-1$
and $|Y|$
are even.

Fill the remaining
positions in
each class with
complete pairs.

The number needed
is precisely
\[
\frac{|X|-1}{2}
+\frac{|Y|}{2}
=\frac{n-1}{2}.
\]

The resulting
weighting is
injective.
Both class
sums are zero,
and its only
zero weight is
at the internal
vertex $r$.

\medskip
\noindent
\textbf{(iii) Odd order,
$k\geq n+2$.}

Let $A$ denote
the odd-sized
class and $B$
the even-sized
class.

When $k$ is odd,
assign $-1$ to
one vertex of
$A$, and assign
$1,-2$ to two
vertices of $B$.

The partial sums
of both classes
equal $-1$.

Avoid the
two pair classes
$\{\pm1\}$
and $\{\pm2\}$.

The remaining
$(n-3)/2$
pairs can be
chosen because
\[
P(k)-2
=\frac{k-5}{2}
\geq\frac{n-3}{2}.
\]

When $k$ is even,
write $k=2q$.
Since $n$ is odd
and $k\geq n+2$,
we have
\[
k\geq n+3,\qquad
q\geq3.
\]

Assign weight
$q$ to one
vertex of $A$,
and weights
$1,q-1$ to
two vertices
of $B$.

Both partial
class sums
equal $q$.

The three
residues are
distinct and
nonzero.

The residue
$q$ is
self-inverse.

Avoid the
two ordinary
pair classes
containing
$1$ and
$q-1$.

There are
enough remaining
pairs because
\[
P(k)-2
=\frac{k-6}{2}
\geq\frac{n-3}{2}.
\]

All constructions
produce distinct
vertex weights
with equal
bipartition sums.
\end{proof}

The exceptional
even-order stars
are handled
directly.

\begin{lemma}[Even-order stars]
\label{lem:stars}
Let $T=K_{1,n-1}$
with even $n\geq4$.

Then $T$ is
$\Z_k$-antimagic
for every
$k\geq n+2$.
\end{lemma}

\begin{proof}
Let $c$ be
the center.

We construct
$n-1$ distinct
nonzero leaf-edge
labels summing
to zero.
Then the center
weight is zero
and differs
from every
leaf weight.

If $k$ is odd,
then
$k\geq n+3\geq7$.

Use the zero-sum
triple
\[
1,\quad2,\quad-3
\]
and complete
the remaining
$n-4$ leaf
positions with
ordinary pairs,
avoiding the
three classes
containing
$1,2,3$.

This is possible
since
\[
P(k)-3
=\frac{k-7}{2}
\geq\frac{n-4}{2}.
\]

If $k$ is even,
write $k=2q$.
Here $q\geq3$.

Use
\[
1,\quad q-1,\quad q,
\]
whose sum is
$2q=0$ modulo $k$.

Complete the
remaining
positions with
ordinary pairs,
avoiding the
two classes
containing
$1$ and $q-1$.

The residue
$q$ belongs to
no ordinary
pair, and
\[
P(k)-2
=\frac{k-6}{2}
\geq\frac{n-4}{2}.
\]

Thus the chosen
leaf-edge labels
are nonzero,
distinct, and
sum to zero.

The center has
weight zero
and all leaf
weights are
distinct and
nonzero, giving
the required
labeling.
\end{proof}

\begin{proof}[Proof of
Theorem~\ref{thm:trees}]
Let $T$ be a tree
of order $n\geq3$.

If $n$ is odd,
Lemma~\ref{lem:weights}
supplies a feasible
injective weighting
for every
$k\geq n+1$.

At $k=n+1$
there is exactly
one zero weight,
at an internal
vertex. For
$k\geq n+2$
all weights
are nonzero.

In both cases
Theorem~\ref{thm:elimination}
gives a nowhere-zero
antimagic labeling.

If $n$ is even
and $T$ is not
a star, apply
Lemma~\ref{lem:weights}(i)
and
Theorem~\ref{thm:elimination}
for every
$k\geq n+1$.

For even-order
stars, apply
Lemma~\ref{lem:stars}
for every
$k\geq n+2$.

These cases
establish the
claimed sufficient
intervals.

Sharpness is
proved below
in
Propositions
\ref{prop:starobstruction}
and \ref{prop:parity}.
\end{proof}

\section{Connected graphs and complete spectra}

We first give
a simple spanning-tree
observation.

\begin{lemma}[Nonstar spanning trees]
\label{lem:spanning}
Let $G$ be a
connected simple
graph of order
$n\geq4$.

If $G$ is not
a star, then
$G$ has a
spanning tree
that is not
a star.
\end{lemma}

\begin{proof}
Choose any
spanning tree
$T$ of $G$.

If $T$ is not
a star, there
is nothing
to prove.

Otherwise let
$c$ be the
center of $T$.

Since $G$
itself is not
a star, there
is an edge
$uv$ between
two leaves of
$T$.

Delete $cu$
from $T$
and add $uv$.

The resulting
graph $T'$
is connected
and has $n-1$
edges, so it
is a spanning
tree.

The vertex $c$
now has degree
$n-2<n-1$.
Every other
vertex has
degree at most
two, and
$n\geq4$.

Thus no vertex
has degree
$n-1$ and
$T'$ is not
a star.
\end{proof}

\begin{proposition}[Large moduli
for connected graphs]
\label{prop:connectedlarge}
Let $G$ be a
connected simple
graph of order
$n\geq3$.

If $G$ is not
an even-order
star, then
\[
[n+1,\infty)\cap\Z
\subseteq\IAM(G).
\]

If $G$ is an
even-order star,
then
\[
[n+2,\infty)\cap\Z
\subseteq\IAM(G).
\]
\end{proposition}

\begin{proof}
For any tree,
both assertions
follow from
Theorem~\ref{thm:trees}.

Suppose $G$
is not a tree.

If $n\geq4$,
then $G$ cannot
be a star.
By
Lemma~\ref{lem:spanning},
choose a nonstar
spanning tree
$T\subseteq G$.

Theorem~\ref{thm:trees}
gives
\[
[n+1,\infty)\cap\Z
\subseteq\IAM(T).
\]

Add the edges
of $E(G)\setminus E(T)$
one at a time.

At every stage
the intermediate
graph is connected
and has order
$n\geq4$.

Thus
Theorem~\ref{thm:knownstability}
preserves every
admissible modulus.

Consequently,
\[
\IAM(T)\subseteq\IAM(G),
\]
and the desired
interval is
contained in
$\IAM(G)$.

It remains to
consider the
non-tree graph
of order three,
namely $K_3$.

For every
$k\geq4$,
label the three
edges of $K_3$
with $1,2,3$.

The three
vertex weights
are
\[
3,\quad4,\quad5
\pmod k.
\]

These residues
are distinct
for every
$k\geq4$.

Hence
\[
[4,\infty)\cap\Z
\subseteq\IAM(K_3).
\]

The order-three
case is therefore
also covered.
\end{proof}

We can now prove
the complete-spectrum
classification.

\begin{proof}[Proof of
Theorem~\ref{thm:complete}]
Let $G$ be connected
of order $n\geq3$.

By the
Jones--Zhang
classification,
Theorem~\ref{thm:knownminimum},
the least
admissible
modulus is
$\mu(G)$, given
by \eqref{eq:minclassification}.

If $G$ is an
even-order star,
then
\[
\mu(G)=n+2.
\]
Proposition~
\ref{prop:connectedlarge}
supplies every
integer $k\geq n+2$.

If $G\cong K_3$,
then
\[
\mu(G)=4,
\]
and the explicit
labeling in the
proof of
Proposition~
\ref{prop:connectedlarge}
supplies every
$k\geq4$.

Suppose next
that $G$ is
not an even-order
star and
\[
n\equiv2\pmod4.
\]

Then
\[
\mu(G)=n+1
\]
by
Theorem~\ref{thm:knownminimum},
and Proposition~
\ref{prop:connectedlarge}
supplies every
$k\geq n+1$.

Finally, in
every other case,
\[
\mu(G)=n.
\]

The modulus
$n$ is admissible
by the published
minimum-modulus
classification.

Proposition~
\ref{prop:connectedlarge}
supplies every
$k\geq n+1$.

Thus in all
cases every
integer at least
$\mu(G)$ is
admissible.

No smaller
integer is
admissible by
the definition
of $\mu(G)$.

Therefore
\[
\IAM(G)=
[\mu(G),\infty)\cap\Z.
\]
\end{proof}

In particular,
the integer-antimagic
spectrum of every
connected simple
graph of order at
least three is
an interval of
integers.

This settles the
interval-spectrum
conjecture stated
in \cite{ORL},
subject to the
published input
theorems and the
new tree proof
given above.

\subsection{Sharpness of the tree bounds}

Although sharpness
of the connected-graph
classification also
follows from
Jones--Zhang, it is
useful to record
two short direct
obstructions.

\begin{proposition}[Even-order stars]
\label{prop:starobstruction}
For every even
$n\geq4$,
\[
n,n+1\notin
\IAM(K_{1,n-1}).
\]
\end{proposition}

\begin{proof}
Let $c$ be
the center.

Modulo $n$,
the $n-1$
leaf labels
must exhaust
all nonzero
residues.

Their sum is
\[
1+\cdots+(n-1)
\equiv\frac n2
\pmod n.
\]

The center
therefore has
the same weight
as the leaf
labeled $n/2$.

Modulo $n+1$,
the modulus
is odd.

The $n-1$
leaf labels
must be all
nonzero residues
except one,
say $a$.

Since the sum
of all nonzero
residues modulo
an odd integer
is zero, the
center weight
is $-a$.

But $-a\neq a$,
so $-a$ is among
the selected
leaf labels.

Again the center
weight collides
with a leaf
weight.
\end{proof}

\begin{proposition}[Parity obstruction]
\label{prop:parity}
If $n\equiv2\pmod4$,
no tree of order
$n$ is
$\Z_n$-antimagic.
\end{proposition}

\begin{proof}
If such a
labeling existed,
its $n$ distinct
vertex weights
would exhaust
$\Z_n$.

Their total
sum would be
\[
0+1+\cdots+(n-1)
\equiv\frac n2
\pmod n.
\]

But for the
bipartition
$V=X\sqcup Y$,
feasibility gives
\[
\sum_{v\in V}
\wt_f(v)
=
2\sum_{x\in X}
\wt_f(x)
\pmod n.
\]

Hence
\[
2\sum_{x\in X}
\wt_f(x)
\equiv\frac n2
\pmod n.
\]

The left
side is even.
Since
$n\equiv2\pmod4$,
the right side
is odd, while
the modulus
$n$ is even.

This congruence
is impossible.
\end{proof}

For $n\geq6$
with $n\equiv2
\pmod4$, the
nonstar path
$P_n$ therefore
fails at
$k=n$.

Together with
Theorem~
\ref{thm:trees},
this shows that
the sufficient
bound $n+1$
cannot be
lowered uniformly
for nonstar trees.

Even-order stars
similarly show
that the universal
bound $n+2$
cannot be
lowered.

\section{Examples and concluding remarks}

We conclude with
three examples
illustrating the
two elimination
mechanisms and
the exceptional
matching configuration.

\begin{example}[Ordinary zero-edge transport]
Let $P_7$ have
vertices
$v_1,\ldots,v_7$
in path order,
and work
modulo $8$.

Prescribe the
vertex weights
\[
w=(4,2,6,1,3,7,5).
\]

The bipartition
sums are
\[
4+6+3+5
\equiv2\pmod8
\]
and
\[
2+1+7
\equiv2\pmod8.
\]

Thus the
weighting is
feasible and
injective.

The inducing
edge labels are
\[
f=(4,6,0,1,2,5).
\]

Swap the weights
at $v_3$ and
$v_5$.

The new vertex
weights are
\[
w'=(4,2,3,1,6,7,5),
\]
and the edge
labels become
\[
f'=(4,6,5,4,2,5).
\]

All edge labels
are now nonzero,
while the vertex
weights remain
distinct.
\end{example}

\begin{example}[Repairing a zero-weight root]
Consider the
path $P_9$
modulo $10$,
with vertices
$v_1,\ldots,v_9$
in order.

Start with
vertex weights
\[
w=(1,3,9,7,0,4,2,6,8).
\]

Both bipartition
sums are zero
modulo $10$.

The induced
edge labels are
\[
f=(1,2,7,0,0,4,8,8).
\]

The unique
zero-weight
vertex is
$r=v_5$,
and both
incident edges
have label zero.

Swap the weights
$7$ and $4$
at $v_4$ and
$v_6$.

The new vertex
weights are
\[
w'=(1,3,9,4,0,7,2,6,8),
\]
and the new
edge labels are
\[
f'=(1,2,7,7,3,4,8,8).
\]

Every edge label
is nonzero.

Direct computation
of the vertex sums
recovers exactly
the distinct weights
in $w'$.

Thus this is a
$\Z_{10}$-antimagic
labeling of $P_9$.
\end{example}

\begin{example}[The exceptional Hall case]
\label{ex:hall}
We exhibit a
configuration in
which the matching
exception genuinely
occurs.

Let $T$ have a
root $r$,
three neighbors
$v_1,v_2,v_3$,
and two leaves
attached to each
neighbor.

Work modulo $17$.

Assign
\[
w(r)=0,
\]
and neighbor
weights
\[
(a_1,a_2,a_3)
=(3,1,5).
\]

Assign leaf
weights
\[
(4,16)
\]
in the first
branch,
\[
(6,14)
\]
in the second,
and
\[
(7,13)
\]
in the third.

All ten
vertex weights
are distinct.

The sum of
the leaf weights
is $60\equiv9$
modulo $17$,
equal to the
sum of the
neighbor weights
$3+1+5=9$.

Thus the
weighting is
feasible.

Since the
two leaves in
each branch
have sum $3$
modulo $17$,
the root-edge
labels are
\[
(\ell_1,\ell_2,\ell_3)
=(0,15,2).
\]

Consequently,
\[
(c_1,c_2,c_3)
=
(\ell_i-a_i)_{i=1}^3
=
(14,14,14).
\]

All positions
forbid the
available
weight $3$.

Hence no
permutation
of the neighbor
weights alone
can make all
root edges
nonzero.

Now swap the
leaf weights
$4$ and $6$
in the first
two branches.

This changes
the root-edge
labels to
\[
(15,0,2),
\]
and the constants
to
\[
(12,16,14).
\]

They are no
longer all
equal.

Assign the
neighbor weights
in the order
\[
(3,5,1).
\]

The resulting
root-edge labels
are
\[
(15,4,15).
\]

All are nonzero
and sum to
zero modulo $17$.

The six leaf-edge
labels remain
the nonzero
residues
\[
6,16,4,14,7,13.
\]

Every edge label
is therefore
nonzero, and
the ten vertex
weights remain
distinct.

This example
shows that the
exceptional
matching case
is real and
that the
cross-branch
swap resolves it.
\end{example}

\subsection*{Constructivity}

The proof gives
a deterministic
polynomial-time
algorithm.

For a tree,
compute its
bipartition
and assign
the feasible
weights of
Lemma~\ref{lem:weights}.

Recover the
unique inducing
edge labeling
by leaf
elimination.

If no vertex
weight is zero,
transport zero
edges until
they disappear.

If one internal
vertex has
weight zero,
first repair
its incident
edges using
the matching
lemma, and
then transport
the remaining
zeros outward.

Each completed
transport strictly
decreases the
number of zero
edges.

A transported
zero moves
along a simple
path and
therefore
requires
$O(n)$ swaps.

There are
at most $n-1$
zero edges,
so $O(n^2)$
transport swaps
suffice.

The matching
step is also
polynomial-time.

The construction
uses only addition
and subtraction
modulo $k$,
so it applies
to composite
moduli as well
as prime ones.

\subsection*{Conclusion}

The principal
new result is
a sharp,
constructive
tree-labeling
theorem.

Its essential
mechanism is
the conversion
of injective
feasible vertex
weights into
nowhere-zero
edge labels
by switching.

The single-zero
case requires
a matching
argument at
the zero-weight
vertex, followed
by outward
zero-edge
transport.

Combined with
the minimum-modulus
classification
of Jones and
Zhang and the
edge-addition
preservation
theorem of
Odabasi, Roberts,
and Low, this
gives a complete
description of
the integer-antimagic
spectra of all
connected simple
graphs of order
at least three.

In particular,
every such
spectrum is an
interval of
integers beginning
at its least
admissible modulus.


\end{document}